\documentclass[twocolumn,11pt]{article}
\usepackage{times}
\newtheorem{definition}{Definition}
\newtheorem{remark}[definition]{Remark}

\newtheorem{theorem}[definition]{Theorem}
\newtheorem{proposition}[definition]{Proposition}
\newtheorem{corollary}[definition]{Corollary}
\newtheorem{example}[definition]{Example}
     \def\({\left(}                  
     \def\){\right)}          
     \def\[{\left[}       
     \def\]{\right]}

     \def\<{\langle}                 \def\wh{\widehat}
     \def\>{\rangle}                 
                 \def\sbs{\subset}

      \def\sbs{\subset}

\newcommand{\eop}{\hfill $\sqcap\!\!\!\!\sqcup$} 
\begin{document}
\global\def\refname{{\normalsize \it References:}}
\baselineskip 12.5pt
%
%
%
\title{\LARGE \bf Geometrical Characterization of RN-operators between Locally Convex Vector Spaces}

\date{}

\author{\hspace*{-10pt}
\begin{minipage}[t]{2.7in} \normalsize \baselineskip 12.5pt
\centerline{OLEG REINOV}
\centerline{St. Petersburg State University}
\centerline{Dept. of Mathematics and Mechanics}
\centerline{Universitetskii pr. 28, 198504 St, Petersburg}
\centerline{RUSSIA}
\centerline{orein51@mail.ru}
\end{minipage} \kern 0in
\begin{minipage}[t]{2.7in} \normalsize \baselineskip 12.5pt
\centerline{ASFAND FAHAD}
\centerline{Government College University}
\centerline{Abdus Salam School of Math. Sciences}
\centerline{68-B, New Muslim Town,  54600 Lahore}
\centerline{PAKISTAN}
\centerline{asfandfahad1@yahoo.com}
\end{minipage}
%
%
\\ \\ \hspace*{-10pt}
\begin{minipage}[b]{6.9in} \normalsize
\baselineskip 12.5pt {\it Abstract:}
  For locally convex vector spaces (l.c.v.s.) $E$ and $F$ and for linear and continuous operator
   $T: E \rightarrow F$ and for an absolutely convex neighbor\-hood $V$ of zero in $F$,
a bounded subset $B$ of $E$ is said to be {\it $T$-V-dentable}\
         (respectively, {\it $T$-V-s-dentable},\, respectively, {\it $T$-V-f-dentable}\,)
 if for any $\epsilon>0$ there exists an $x\in B$ so that
 $$
 x\notin \overline{co}\, (B\setminus T^{-1}(T(x)+\epsilon V))
 $$
 (respectively, so that
  $
 x\notin s$-$co\, (B\setminus T^{-1}(T(x)+\epsilon V)),$\, 
respectively, so that
  $
 x\notin {co}\, (B\setminus T^{-1}(T(x)+\epsilon V))\, ).
 $
Moreover, $B$ is called {\it $T$-dentable}\ (respectively, {\it $T$-s-dentable}, {\it $T$-f-dentable}\,) 
if it is {\it $T$-V-dentable}\ (respectively, {\it $T$-V-s-dentable}, {\it $T$-V-f-dentable}\,) for every absolutely convex neighbor\-hood $V$ of zero in $F.$
 RN-operators between locally convex vector spaces have been introduced in [5]. 
  We present a theorem 
   which says that, for a large class of l.c.v.s. $E, F,$ if   
   $T: E \rightarrow F$ is a linear continuous map, then the following are equivalent:\\
         1) \  $T \in RN(E,F).$\\
         2)\ Each bounded set in $E$ is $T$-dentable.\\
         3)\ Each bounded set in $E$ is $T$-s-dentable.\\
         4) \ Each bounded set in $E$ is $T$-$f$-dentable.\\                                               
Therefore, we have a generalization of Theorem 1 in [8], 
 which gave a geometric characterization of RN-operators between Banach spaces.
\\ [4mm] {\it Key--Words:}
dentability,  dentable set, locally convex space, Radon-Nikodym operator
\end{minipage}
\vspace{-10pt}}

\maketitle

\thispagestyle{empty} \pagestyle{empty}
%
%
\section{Introduction}
\label{S1} \vspace{-4pt}

 We use standard notations and terminology from the theory of operators in Banach spaces
 and measure theory (see, e.g., \cite{Diestel},\cite{Dunford}).
 Our main references on the theory of the locally convex vector spaces is \cite{H. H. Schaefer}.   
 Banach spaces are usually denoted by the letters $X,Y, \dots,$ and the usual notations for the l.c.v.s.
 (locally convex vector spaces) are $E, F, \dots.$ All the spaces we consider in the paper are over the field of the real numbers.

 Some of the main notions and notations we will use are ones of the convex, s-convex and closed convex hulls of sets
 in Banach or l.c.v. spaces. Namely,
 
 $co(A)$ (respectively, $\overline{co}(A),$ respectively, $s$-$co(A))$ denotes the convex hull
 (respectively, the closure of the convex hull,  the s-closure of the convex hull) of the set $A$.
 The s-closure is defined as
 $ s$-$co(A):= \bigl\{\sum\limits_{i\ge 1} a_{i}x_{i}, $
where $ x_{i}\in A, $  $ a_{i}\ge0 $ and  $ \sum\limits_{i\ge 1} a_{i}=1\bigr\}. $

 Let us recall some more notions:

\begin{definition}\label{1.1}
Let $X$ be a Banach space and $(\Omega,\Sigma)$ a measureable space, consider $m: \Sigma \to X$.
$m$ is called a vector measure if for every sequence $\{A_i\}_{1}^{\infty}$ of pairwise disjoint sets
from $\Sigma$ one has $m(\bigcup_{i=1}^{\infty}A_{i})=\sum_{i=1}^{\infty}m(A_{i})$.
\end{definition}
\begin{definition}
Let $m: \Sigma \to X$ be a vector measure. Variation of a vector measure is a non-negative,
extended real valued function with value on set $A\in \Sigma$ is
$$
\ |m|(A)=\sup\limits_{\prod}\sum\limits_{A_i\in\prod}\|m(A_i)\|
$$
where $\prod$ denotes all finite partitions of $A$ with pairwise disjoint sets in $\Sigma$; 
$m$ is called \textit{vector measure of bounded variation} if $|m|(\Omega)$ is finite.
\end{definition}
 
Let $(\Omega,\Sigma, \mu)$ be a measure space.

\begin{definition}
A function $\ f:\Omega \to X $ is \textit{$\mu$-measureable} if there exists a sequence 
of simple functions $\{f_n\}_{n=1}^{\infty}$ with $\lim \limits_{n} \|f_n-f \|=0$  \,\,$\mu$-a.e. {\rm(see \cite{Diestel},page 41)}.
\end{definition}
\begin{definition}
 A $\mu$-measureable function
$
f:\Omega \to X
$
 is called \textit{Bochner integrable} if there exists a sequence of simple functions $f_n$  such that
 $ \lim \limits_{n}\int_\Omega \|f_n-f\| d\mu=0$ {\rm(see \cite{Diestel},page 44)}; in this case
 we set $\int_{A}f=\lim \limits_{n}\int_{A}f_n$ for every $A\in\Sigma$.
  \end{definition}

 The following notions have been introduced in \cite{RN} and, independently, in \cite{WLinde}.
 \begin{definition}
A linear continuous operator $T:X \rightarrow Y$  is said to be RN-operator $(T \in RN(X,Y)$),
 if for every $X$-valued measure $m$ of bounded variation,
 for any measure space $(\Omega,\Sigma, \mu)$ with 
   $m<<\mu$
 there exists a function $f:\Omega$ $\rightarrow Y $ which is Bochner integrable and such that
 $$T(m)(A)=\int_{A}f d\mu$$
  for every $A\in \Sigma$.
 \end{definition}

\begin{definition} Let $T:X \rightarrow Y$ be a linear continuous operator A bounded subset $B$
of Banach space $X$ is called $T$-dentable (respectively, $T$-f-dentable, $T$-s-dentable),
if for every $\epsilon>0$
there exists $x\in B$ such that
$$
   x\notin \overline{co}\, (B\setminus T^{-1}(D_{\epsilon}(T(x)))).
   $$
    (respectively, so that
  $$
 x\notin {co}\,  (B\setminus T^{-1}(D_{\epsilon}(T(x))\, )),
 $$
 resp.,
   $$
 x\notin {s-co}\,  (B\setminus T^{-1}(D_{\epsilon}(T(x))\, ))),
 $$
\end{definition}

The geometrical characterization of RN-operators between Banach spaces has been given by following theorem in \cite{Reinov Geometric}.
\begin{theorem}\label{geometry}
Let $X$, $Y$ be Banach spaces and $T:X \rightarrow Y$, linear and continuous then the following statements are equivalent:
                    
1,\  $T \in RN(X,Y)$

2.\  Each bounded set in $X$ is $T$-dentable.

3.\ Each bounded set in $X$ is $T$-$s$-dentable.

4.\ Each bounded set in $X$ is $T$-$f$-dentable.
\end{theorem}

  Our aim is to generalize the above theorem for operators between some l.c.v. spaces.
  Firstly, we shall define the notion of RN-operators between locally convex vector spaces.
  For this, we recall the definitions of Banach spaces of type $E_B$ and $\hat{E_{V}}.$
 Then we  define
  the notions of $T$-dentability, $T$-s-dentability and $T$-f-dentability 
 for operators in l.c.v.s. 
(in the next section).

 Let $E$ be a l.c.v.s. and $B\subseteq E$ be bounded and absolutely convex. 
 Let $E_{B}=\bigcup_{n=1}^{\infty}nB$. Define  Minkowski function on $E_{B}$ w.r.t. $B:$

   $${\|x\|}_{\rho_{B}}={\rho} _{B}(x)=\inf\{\lambda \geq 0: x\in \lambda B \}.$$

It is a seminorm in general but we show that when $B$ is bounded and absolutely
convex then it is a norm. Suppose $x\in E_{B}$ such that
$\rho_{B}(x)= 0$ so $\inf\{\lambda \geq 0: x\in \lambda B \}=0$  this implies
for any $\epsilon>0$  $\exists$ $\lambda \geq 0$ such that $x\in \lambda B$
and $\lambda\le \epsilon$, take $\epsilon=\frac{1}{n} $, so
$x\in \frac{1}{n}B$ $\forall$ $n\in \mathrm{N}$ this implies
$x\in\bigcap_{n=1}^{\infty}\frac{1}{n}B=\{0\}$. So $x=0$ and $(E_{B},\rho_{B})$ is a
normed space and its completion $\hat{E_{B}}$ is Banach space in $E$.
Moreover, if $B$ is complete then $E_B$ is Banach space.

Similarly for  an absolutely convex open neighbor\-hood of zero $V=V(0)$ we set
 $E_{V}=\bigcup_{n=1}^{\infty}nV=E.$ Similarly we define Minkowski function on $E_{V}$ w.r.t $V$, so that for $x \in E $ the semi norm of $x$ is


$${\|x\|}_{\rho_{V}}= \rho_{V}(x)= \inf\{\lambda \geq 0: x\in \lambda V \}$$.\\

We identify  two elements $x,y\in E$ w.r.t. $\rho_{V},$
obtaining a quotient $E_V$ with corresponding elements $\hat x, \hat y,$ by
$\hat{x}=\hat{y}$  iff  $\rho_{V}(x-y)=0$. 
We get a normed  space $\hat{E_{V}}$ and its completion gives us a Banach space $(\hat{E_{V}},\rho_{V}).$\\
 We will use the following important theorem (see \cite{a and r}) in the proof of our main theorem.
\begin{theorem} {\rm(see \cite{a and r})} \label{pro}
Let $E$ be a locally convex  vector space, $V=V(0)$ be an absolute convex neighbor\-hood of $0$,
 let $ B \subseteq E $ be a closed, bounded, convex, sequentially complete and metrizable subset, 
 The following are equivalent:\\
  (i) $B$ is subset $V$-dentable.\\
 (ii) $B$ is subset $V$-f-dentable.\\
\end{theorem}

It has been shown in  \cite{a and r} that it follows from the above theorem:

{\it
 Let $E$ be a locally convex  vector space and let $ B \subseteq E $ have the following properties:\\
 (i)  $B$ is closed, bounded, convex and sequentially complete,\\
(ii) for every bounded $ M \subseteq E $ and for $ x\in \overline{M}$ there exist a sequence
${x_{n}}\in M$ such that $\lim \limits_{n} {x_{n}}=x,$ \\
(iii) each separable subset of $B$ is metrizable.\\
 Then the  following are equivalent:\\
  (i) $B$ is subset $V$-dentable,\\
 (ii)$B$ is subset $V$-s-dentable,\\
 (ii)$B$ is subset $V$-f-dentable.\\
}

We will say that {\it
a locally convex  vector space $E$ is an SBM-space if

$(i)$\
every closed bounded convex subset  of $E$
is sequentially complete,

$(ii)$\
for every bounded $ M \subseteq E $ and for $ x\in \overline{M}$ there exists a sequence
${x_{n}}\in M$ such that $\lim \limits_{n} {x_{n}}=x$  and

$(iii)$\
 each separable bounded subset of $E$ is metrizable.
}

Therefore, if $E$ is an SBM-space, then every bounded subset of $E$ is $V$-dentable iff
every bounded subset of $E$ is $V$-s-dentable iff
every bounded subset of $E$ is $V$-f-dentable.

All quasi-complete (BM)-spaces \cite{Saab}, in particular, all Fr\'{e}chet spaces are SBM-spaces.   


\section{Main Results}
\vspace{-4pt}

\begin{definition}
         Let $E$ and $F$ be locally convex vector spaces (l.c.v.s.), let $T: E \rightarrow F$
         be a linear and continuous operator and let $V$ be an absolutely convex neighbor\-hood
         of zero in $F.$ A bounded subset $B$ of $E$ is said to be {\it $T$-V-dentable}\
         (respectively, {\it $T$-V-s-dentable},\, respectively, {\it $T$-V-f-dentable}\,)
 if for any $\epsilon>0$ there exists an $x\in B$ so that
 $$
 x\notin \overline{co}\, (B\setminus T^{-1}(T(x)+\epsilon V))
 $$
 (respectively, so that
  $$
 x\notin {s-co}\, (B\setminus T^{-1}(T(x)+\epsilon V))\, ,
 $$
respectively, so that
  $$
 x\notin {co}\, (B\setminus T^{-1}(T(x)+\epsilon V))\, ).
 $$
\end{definition}

\begin{definition}
Let $T: E \rightarrow F$ be a linear and continuous operator. A bounded subset $B$ of $E$
is said to be {\it $T$-dentable}\ (respectively, {\it  $T$-s-dentable,\, $T$-f-dentable}\,) if for every absolutely
convex neighbor\-hood $V$ of zero in $F$ there exists an $x\in B$ such that
 $$
 x\notin \overline{co}\, (B\setminus T^{-1}(T(x)+ V))
 $$
 (respectively, so that
  $$
 x\notin {s-co}\, (B\setminus T^{-1}(T(x)+ V))\, 
 $$
and, respectively, so that
  $$
 x\notin {co}\, (B\setminus T^{-1}(T(x)+ V))\, ).
 $$
\end{definition}

\begin{remark}
 From the above two definitions it is clear that $B \subseteq E$ is $T$-dentable if and only if
 for every $V$ it is $T$-V-dentable.
 The same is true for corresponding properties of s-dentability and f-dentability.
\end{remark}

The following is the definition of an RN-operator between locally convex vector spaces,
and this is our main definition.

\begin{definition} {\rm(see \cite{O n A pre})}
Let $T:E \rightarrow F$ be linear and continuous (in l.c.v.s.).
$T \in RN(E,F)$ (a Radon-Nikodym operator or RN-operator) if for every complete,
absolutely convex and bounded set $B\subseteq E$ and for any absolutely convex
neighbor\-hood $V \subseteq F$ of zero the natural operator
$
\Psi_{V}\circ T \circ \phi_{B}:\hat{E_{B}} \rightarrow E \rightarrow F \rightarrow \hat{F_{V}}$
is $RN$-operator between Banach spaces $\hat{E_{B}}$ and $\hat{F_{V}}.$
\end{definition}

Here, $\Psi_V$ and $\phi_{B}$ are the natural  maps (cf. definitions in Section 1).
\vskip 0.1cm

              \noindent
\begin{remark}                  {
For the operators in Banach spaces our definition of $RN$-operators coincides with the original definitions from  {\rm\cite{RN},
\cite{WLinde}}. This must be clear.
            }
     \end{remark}
\vskip 0.1cm

              \noindent
\begin{remark}                {
The usual definition of a weakly compact operator between locally convex spaces is:
$T: E\to F$ is weakly compact if $T$ maps a neighborhood of zero in $E$ to
a relatively weakly compact subset of $F.$
The analogous definition can be given for the compact operators. For the weakly compact case this means that

(1)\
There exists an absolutely convex neighborhood $V=V(0)$ in $E$ such that if $D:=T(V),$
then $D$ is bounded in $F$ and  the natural injection $\phi_D: F_D\to F$ is weakly compact.

Since every weakly compact operator in Banach spaces is Radon-Nikodym,
we get from the definition that every weakly compact operator in l.c.v.s. is an RN-operator.

Following this way,  we can defined also a class of "bounded RN-operators". Namely, let $T$ maps $E$ to $F.$ $T$ is
a bounded RN-operator if  $T$ takes a neighborhood $V$ to a bounded subset of $F$ and
the natural map $\phi_{T(V)}: F_{T(V)}\to F$ is "Radon-Nikodym".
But how to understand "Radon-Nikodym" in this case where an operator maps a normed (or Banach) space to a locally convex space?
We can go by a geometrical way (saying that the image of $\phi_{T(V)}$ is subset s-dentable. Here the map $\phi_{T(V)}$
is one-to-one, and we can follow the assertion from \cite{Reinov Geometric} for operators in Banach spaces: if $U:X\to Y$ is one-to-one then           
$U$ is RN iff the U-image of the unit ball is subset s-dentable).
Another way is just to apply the main definition from this paper.
Or, we can go by a traditional way: for $T: X\to F$ with $X$ Banach and $F$ locally convex,
say that $T$ is of type RN if for every operator $U$ from an $L_1$-space to $X$
the composition $TU$ admits an integral representation with a (strongly) integrable $F$-valued
function. But in this case we are to give a good definition  of the "integrability"
of an $F$-valued function. All these are the topic for the further considerations in another paper.

Let us mention that in every "right" definition of RN-operator there must be an "ideal property":
$ATB$ is RN for all linear continuous $A,B$ if $T$ is RN. Thus, if we apply the definition 12  
above as the main definition (in this paper), then every operator "of type RN" considered
above in this remark is "right" Radon-Nikodym.
This is one of the reason that here we will deal only with Definition 12.     
            }
            \end{remark}
\vskip 0.1cm

If $E=X$ is a Banach space then $T\in RN(X,F) $ iff for every (a.c.) neighborhood $V=V(0)\sbs F$ the composition $\Psi_VT$
belongs to $RN(X,\wh{F}_V).$
If $F$ is a Banach space too, then, evidently, $T\in RN$ in the sense of Definition 12    
iff it is an RN-operator in the usual sense of \cite{RN} and \cite{WLinde}.    
If $E$ is a l.c.v.s. and $F=Y$ is Banach, then $T\in RN(E,Y)$ iff for every bounded complete absolutely convex subset
$B\sbs E$ the natural map $T\phi_B: E_B\to E\to Y$ is Radon-Nikodym, what means that $T\in RN(E,Y)$
iff for any finite measure space $(\Omega, \Sigma,\mu)$ and for every $\mu$-continuous $E$-valued measure
$\bar m:\Sigma \to E$ with bounded $\mu$-average $\{\frac{\bar m(A)}{\mu(A)}:\ A\in \Sigma,\, \mu(A)\neq0\}$
the measure $T\bar m:\Sigma\to Y$ has a Bochner derivative with respect to $\mu.$
From this it follows that an operator $T$ between locally convex spaces $E,F$ is Radon-Nikodym (in the sense of Definition 12)  
iff
for any finite measure space $(\Omega, \Sigma,\mu)$ and for every $\mu$-continuous $E$-valued measure
$\bar m:\Sigma \to E$ with bounded $\mu$-average, for every a.c. neighborhood $V=V(0)\sbs F$
the measure $\Psi_VT\bar m:\Sigma\to \wh{F}_V$ has a Bochner derivative with respect to $\mu.$
\vskip 0.1cm

\begin{proposition}\label{prop1}
Let $\Psi_{V}:E \rightarrow \hat{E_{V}}$, where $V=V(0)$ is an absolutely convex open
neighbor\-hood of zero then
${\Psi_{V}}^{-1}(D_{1}(0))=V$, where $D_{1}(0)$ is open unit ball in $\hat{E_{V}}.$
\end{proposition}

\noindent
{\bf Proof:} \  
$"\supseteq"$:\
Let $ x\in V$, then by definition of $\|. \|_{\hat{E_{V}}}$,  $\|\Psi_{V}(x)\|< 1 $  so
 $V \subseteq {\Psi_{V}}^{-1}(D_{1}(0))$%
 
$ "\subseteq"$:\
Let $x\in{\Psi_{V}}^{-1}(D_{1}(0));$ then $\Psi_{V}(x)\in D_{1}(0).$ This
implies $\|\Psi_{V}(x)\|=1-c <1$ for some $c>0$, so
$\inf\{\lambda > 0: x\in \lambda V \}=1-c <1 .$ This implies
$x\in {(1-{\frac{c}{2}})}V$. Since $V$ is absolutely convex so
$\lambda V \subseteq V $ for every $|\lambda|<1$, therefore $x\in V$.
Hence ${\Psi_{V}}^{-1}(D_{1}(0))=V.$   \eop

%
%



\begin{theorem}
Let $E$, $F$ be locally convex vector spaces and $T: E \rightarrow F$ be a linear continuous operator. Consider the following conditions:
                     
1)\  $T \in RN(E,F).$

2)\ Each bounded set in $E$ is $T$-dentable.

3)\ Each bounded set in $E$ is $T$-$s$-dentable.

4)\ Each bounded set in $E$ is $T$-$f$-dentable
                                                                           
\noindent                                                                           
We have      $1)\Leftrightarrow 4),$\, $2)\Rightarrow 3)\Rightarrow 4).$
If the space   $E$    is an SBM-space then      all the conditions are equivalent.
\end{theorem}

\noindent
{\bf Proof:} \  
$1)\Rightarrow 4).$\ 
          Let $B_{0}\subseteq E$ be bounded and $V\subseteq F$ be an absolutely convex
     neighbor\-hood of $0$; put $B=\overline{\Gamma(B_{0})}$.
     By assumption in $1),$ the composition operator
     $\Psi_{V}\circ T \circ \phi_{B}:E_{B} \rightarrow E \rightarrow F \rightarrow \hat{F_{V}}$
     is an $RN$-operator from $E_{B}$ to $\hat{F_{V}}.$ By theorem \ref{geometry} each bounded subset
     in $E_{B}$ is $\Psi_{V} T \phi_{B}$-f-dentable. In particular, $B_0$ is
     $\Psi_{V} T \phi_{B}$-f-dentable. Let $\epsilon>0$, there exists $x\in B_0$ such that
   $$
  (i)\, \ x\notin {{co}}\, (B_{0}\setminus (\Psi_{V} T \phi_{B})^{-1}{(D_{\epsilon}(\Psi_{V} T \phi_{B}(x)))}).
   $$


\noindent
  Proposition \ref{prop1} implies that ${\Psi_{V}}^{-1}(D_{\epsilon}(\Psi_{V} T \phi_{B}(x)))={\Psi_V}^{-1}(\Psi_V(T(x)))+\epsilon V$, so we get\\
$$
(ii)\,   \ x\notin {co}\, (B_{0}\setminus  {T}^{-1}{(T(x)+\epsilon V)},
   $$
 otherwise there exist $x_1, x_2$ ..., $x_n$  in $B_0$ where $x_{k}\notin {{T}^{-1}{(T(x)+\epsilon V)}}$ and $\lambda_1, \lambda_2$ ...,$\lambda_n$ where $0\leq \lambda_k \leq 1$ and $\sum_{k=1}^{n} \lambda_k = 1$ such that $\sum_{k=1}^{n} \lambda_k x_k=x $.
Since $x_{k}\notin {{T}^{-1}{(T(x)+\epsilon V)}}$ so $T(x_{k})\notin {(T(x)+\epsilon V)}$ hence  $T(x_{k})\notin {\Psi_V}^{-1}(\Psi_V(T(x)))$ for each
$k= 1,2$ ... $n$. So, we get $x_{k}\notin {\phi_B}^{-1}{T^{-1}\Psi_{V}}^{-1}(D_{\epsilon}(\Psi_{V} T \phi_{B}(x)))$ and $\sum_{k=1}^{n} \lambda_k x_k= x $ 
which is contradicts $(i)$, therefore $(ii)$ holds and hence $B_0$ is $T$-V-f-dentable, which proofs $4)$.

$4)\Rightarrow 1).$\
 Let $B\subseteq E$ be a bounded, absolutely convex and complete subset and $V\subseteq F$ be an absolutely convex neighbor\-hood of $0.$ We need to show that
 the operator $\Psi_{V} \circ T \circ \phi_{B} : E_{B} \rightarrow E \rightarrow F \rightarrow \hat{F_{V}}$ is $RN$.
 By theorem \ref{geometry}, it is enough to show that each bounded $B_0$ is $\Psi_{V} T\phi_{B}$-f-dentable.
 Let $\epsilon>0;$ by assumption in $4)$ there exists an  $x\in B_{0}$ such that
  $$
 x\notin {co}\, (B_{0}\setminus T^{-1}(T(x)+ \epsilon V))\, ),
 $$
 which, by above proposition, gives that
$$
 x\notin {co}\, (B_{0}\setminus (\Psi_{V} T\phi_{B})^{-1}(D_{\epsilon}(\Psi_{V} T \phi_{B}(x)))\,).
 $$
This proofs $1)$.

$2)\Rightarrow 3)\Rightarrow 4).$\
 Follows from definitions.
\smallskip

Suppose now that $E$ is an SBM-space.

$4)\Rightarrow 2).$\
 Let $B\subseteq E$ be closed bounded and convex and let $V$ be an absolutely convex neighbor\-hood of zero in $F$. We show that $B$ is subset $T$-V-dentable.
 Clearly, $B_0\subset B$ is $T$-V-f-dentable if and only if $B_0$ is $U=T^{-1}(V)$-f-dentable.
 From assumption, it follows that $B$ is subset $U$-f-dentable. Applying the consequence of theorem \ref{pro},
  we get that $B$ is subset $U$-dentable, or $B$ is subset $T$-V-dentable, which proofs $2)$.
\eop



\begin{example}
{\rm
Consider any uncountable set $\Gamma$ and the classical Banach spaces $c_0(\Gamma)$ and $l_1(\Gamma);$
$c_0(\Gamma)^*=l_1(\Gamma).$
The closed unit ball $B$ of $l_1(\Gamma)$ is weak$^*$ compact and an Eberlein compact in this
weak$^*$ topology (see, e.g., \cite{DiesGBS}).   
This implies that every separable subset of $(B, w^*)$ is metrizable and that $(B,w^*)$ has the Frechet-Uryson property,
i.e. every point in the closure of a subset is a limit of a sequence of this subset.
Thus, the space $(l_1(\Gamma),w^*)$ is of type SBM.           
The Banach space $l_1(\Gamma)$ has the RNP. Therefore, every bounded subset of the space is dentable (s-dentable, f-dentable),
and thus every bounded subset of $(l_1(\Gamma), w^*)$ is ($w^*)$-dentable too. The identity map
$(l_1(\Gamma), w^*)\to (l_1(\Gamma, w^*))$ is Radon-Nikodym in the sense of Definition 12.    
Moreover the space has the RNP in the sense of the paper \cite{Saab}.
This directly follows by Theorem 16,         
but not {\it directly}\, from the results of \cite{Saab}. On the other case, for $l_1(\Gamma)$ this fact is trivial.
We can also say (by the same considerations) that if $X$ is any weakly generated Banach space (see \cite{DiesGBS}) with $X^*\in RN,$
then the space $(X^*, \sigma(X^*,X))$ has all just mentioned properties (clearly, or by application of Theorem 16).  } 
\end{example}

Finally, if  $X$ is WCG, then $(X^*, \sigma(X^*,X))$ is SBM. So, the theorem can be applied for all such spaces.

On the other hand, there is an example (see \cite{Egge}, Theorem 4), which shows that there exists a separable
Banach space $Y$ such that the space $Y_{\sigma:=}(Y, \sigma (Y,Y^*))$ does not have the RNP (in the sense
of \cite{Saab})
 and every bounded subset of which
(of $Y_{\sigma})$ is $s$-dentable. Thus, the above theorem is not true for general l.c.v.s.

Let us note that a l.c.v.s. $E$ has the RNP in sense of \cite{Saab} iff the identity map $E\to E$ is an RN-operator.
Therefore, we obtain a generalization of a theorem from \cite{Saab}:

\begin{corollary}
For every SBM-space $E$ (in particular, for every quasi complete l.c.v.s. with metrizable bounded subsets, or
for every Frechet space) the following are equivalent:

                                                                           
1.\  $E$ has the RNP of  \cite{Saab}.

2.\ Each bounded set in $E$ is dentable.

3.\ Each bounded set in $E$ is $s$-dentable.

4.\   Each bounded set in $E$ is $f$-dentable.                                                                         
\end{corollary}


\vspace{10pt} \noindent
{\bf Acknowledgements:} \ The research was supported by the Higher  Education Commission of Pakistan
and by the Grant Agency of RFBR (grant No. 15-01-05796).


\end{document}